\documentclass[11pt]{amsart}
\usepackage{latexsym, amssymb, amsmath}
\def\R{\mathbb R}

\def\Z{\mathbb Z}
\def\RP{\mathbb {RP}}

\def\vol{\mathrm{Vol}}

\def\d{\mathrm{div}}
\def\dist{\mathrm{dist}}

\def\as {asymptotically flat }
\def\imc{inverse mean curvature }

\def\l{\mathrm{loc}}
\newtheorem*{conj}{Conjecture}
\newtheorem{thm}{Theorem}[section]
\newtheorem{lemm}[thm]{Lemma}
\newtheorem{cor}[thm]{Corollary}

\theoremstyle{remark}
\newtheorem{rmk}[thm]{Remark}
\theoremstyle{definition}
\newtheorem{defi}[thm]{Definition}

\title{ $3$-manifolds with Yamabe invariant greater than that of $\RP^3$}

\author{Kazuo Akutagawa${}^*$}
\email{akutagawa$_{-}$kazuo@ma.noda.tus.ac.jp}
\address{Department of Mathematics, Tokyo University of Science,
Noda 278-8510, Japan}

\author{Andr\'e Neves${}^{\dagger}$}
\email{aneves@math.princeton.edu}
\address{Department of Mathematics, Princeton University,
Princeton, NJ 08540, USA}

\thanks{${}^*$\
supported in part by the Grants-in-Aid for Scientific Research (C),
Japan Society for the Promotion of Science, No.~16540059.\\
\quad\ ${}^{\dagger}$\ supported in part by
FCT-Portugal, Grant BD/893/2000.}

\date{October, 2004.}
\begin{document}
\maketitle \markboth{ $3$-manifolds with Yamabe invariant greater than that of $\RP^3$}
{Kazuo Akutagawa and Andr\'e Neves}

\begin{abstract}

We show that, for all nonnegative integers $k,\  \ell,\  m$ and $n$,
the {\em Yamabe invariant} of
$$\# k(\RP^3) \# \ell(\RP^2\times S^1)
\# m(S^2 \times S^1) \# n(S^2 \tilde{\times} S^1)$$ is equal to the {\em Yamabe invariant} of $\RP^3$, provided
$k + \ell \geq 1$. We then complete the classification (started by Bray and the second author) of all closed
$3$-manifolds with Yamabe invariant greater than  {that of }$\RP^3$. {More precisely, we show } that such
manifolds are either $S^3$ or finite connected sums $\# m(S^2 \times S^1) \# n(S^2 \tilde{\times}
S^1)$, where $S^2 \tilde{\times} S^1$ is the nonorientable $S^2$-bundle over $S^1$.

A key ingredient {  is}  Aubin's Lemma \cite{Aubin-1}, which says that if the {\em Yamabe constant} is positive,
then it is strictly less than the {\em Yamabe constant} of any of its non-trivial finite conformal coverings.
This lemma, combined with  {  inverse mean curvature flow} and with analysis of the Green's functions for the
conformal Laplacians on specific finite and normal infinite Riemannian coverings, will allow us to construct a
family of nice test functions on the finite coverings and thus prove the desired result.

\end{abstract}
\maketitle

\section{Introduction}
It is well known that, for a closed smooth  {n-}manifold $M $  {with } $n \geq 3$, a Riemannian metric is
Einstein if and only if it is a critical point of the normalized Einstein-Hilbert functional $E$ on the space
$\mathcal{M}(M)$ of all Riemannian metrics on $M$
$$ E : \mathcal{M}(M) \rightarrow {\mathbb R}, \quad
g \mapsto E(g) = \frac{\int_MR_g~dV_g}{\textrm{Vol}_g(M)^{(n-2)/n}}.$$ Here, $R_g, dV_g$ and $\textrm{Vol}_g(M)$
denote respectively the scalar curvature, the volume form of the metric $g$, and the volume of $(M, g)$. { Next,
we give the definition of a natural differential-topological invariant $Y(M)$ of $M$, called the {\it Yamabe
invariant} {or $\sigma$-invariant}}.  {Denote the conformal class of a Riemannian metric $g$ by}
$$C = [g] := \{ \textrm{e}^{2f}\cdot g\ |\ f \in C^{\infty}(M) \}$$ and {  the space of all conformal structures
by }$\mathcal{C}(M)$.
{  The} Yamabe invariant $Y(M)$ is defined by
$$
Y(M) := \sup_{C \in \mathcal{C}(M)} \inf_{g \in C} E(g).
$$
Contrarily to what was stated in \cite{BN}, this invariant was
introduced independently by
O.~Kobayashi \cite{Kobayashi} and Schoen \cite{Schoen-3}.

To put this definition in context, first recall that $E$ is neither bounded from above nor { from} below for any
$M$ and so a minimization method cannot be applied to $E$. Nevertheless, the study of the second variation of
$E$ done in \cite{Koiso, Schoen-2} (cf.\ \cite{Besse}) leads naturally, as we will argue next, to the mini-max
definition of the Yamabe invariant.

Let $\mathcal{M}_1(M)$ denote the unit volume metrics of $M$.
Due to results of Obata \cite{Obata}
and Aubin \cite{Aubin-2}, an Einstein metric $h$
in $\mathcal{M}_1(M)$ globally minimizes $E$ in $\mathcal{M}_1(M) \cap [h]$.
Moreover, if $(M, [h])$ is not conformally equivalent to $(S^n, [g_0])$,
$E$ is strictly globally minimized among variations in $T_h(\mathcal{M}_1(M) \cap [h])$.
Here, $g_0$ is the canonical metric of constant curvature one on the standard $n$-sphere $S^n$.
On the other hand, if we consider the infinite dimensional subspace $V$
consisting of infinitesimal deformations orthogonal to both $T_h[h]$
and the subspace of infinitesimal deformations arising from the diffeomorphism group $\textrm{Diff}(M)$ of $M$,
$E$ is locally maximized among variations in $V$
modulo a finite dimensional subspace (possibly the zero subspace).
These observations roughly imply that critical points are always saddle points
over the moduli space $\mathcal{M}_1(M)/\textrm{Diff}(M)$.

The restriction of $E$ to any given conformal class $[g]$ is always bounded from below.
Hence, we can consider the following conformal invariant
\begin{align*}
& Y(M,[g]) := \inf\{E(\tilde{g})\ |\ \tilde{g} = u^{4/(n-2)}\cdot g \in [g], \
u > 0 \mbox{ in } C^{\infty}(M) \} \\
& = \inf\left\{Q_g(u) :=
\frac{\int_M (\frac{4(n-1)}{n-2}|\nabla u|^2 + R_g u^2)dV_g}
{\left(\int_M u^{2n/(n-2)}~dV_g\right)^{(n-2)/n}}\
\Bigg| \  u>0 \mbox{ in } C^{\infty}(M)\right\},
\end{align*}
which is called the {\it Yamabe constant} of $[g]$.
A remarkable theorem \cite{Yamabe, Aubin-2, Aubin-Book, Schoen-1}
(cf.~\cite{LP, SY-Book}) of
Yamabe, Trudinger, Aubin, and Schoen asserts that
each conformal class $[g]$ contains metrics $\check{g}$,
called \textit{Yamabe metrics}, which realize the minimum
$$Y(M,[g]) = E(\check{g}).$$
These metrics must have constant scalar curvature
$$R_{\check{g}} = Y(M, [g])\cdot \textrm{Vol}_{\check{g}}(M)^{-2/n}.$$
Aubin \cite{Aubin-2} (cf.~\cite{Aubin-Book}) also proved that
$$Y(M, C) \leq Y(S^n, [g_0]) = n(n-1)\textrm{Vol}_{g_0}(S^n)^{2/n}$$
for any $C \in \mathcal{C}(M)$.
This implies that
$$Y(M) = \sup_{C \in \mathcal{C}(M)} Y(M, C)
\leq Y(S^n) = n(n-1)\textrm{Vol}_{g_0}(S^n)^{2/n}.$$

There is a substantial body of work on Yamabe invariants
and we will briefly review some significant results
(see {  the} surveys \cite{Kazuo, LM, LeBrun-4, Stolz} for example).
The Yamabe invariant $Y(M)$ is positive if and only if
$M$ carries a metric of positive scalar curvature.
Some classification theorems for manifolds with positive Yamabe invariant
have been obtained in \cite{GL-1, GL-2, SY-1, SY-2, Stolz} which show, for instance,
that the $n$-torus $T^n$ has $Y(T^n) = 0$.
There are also some useful surgery theorems
\cite{Kobayashi, Petean, PY} for Yamabe invariants.
LeBrun \cite{LeBrun-3} proved that, for compact complex surfaces,
the Yamabe invariant is positive, zero, or negative,
if and only if the Kodaira dimension is $- \infty, 0$ or $1$, or $2$, respectively.
LeBrun and collaborators \cite{GL, IL, LeBrun-1, LeBrun-2} also computed the Yamabe invariants
for a large class of $4$-manifolds,
including K\"ahler surfaces $X$ with $Y(X) < 0$ and $0 < Y(X) < Y(S^4)$.
In \cite{LeBrun-4}, it is constructed a {  K\"{a}hler} surface which admits Einstein metrics
but the Yamabe invariant is not achieved by any Yamabe metric.

Regarding $3$-manifolds,
the following conjecture by Schoen \cite{Schoen-2}
is one of the most fundamental unsolved problems related to Yamabe invariants.

\begin{conj} If $M$ is a closed $3$-manifold
which admits a constant curvature metric $h$,
then
$$Y(M) = E(h).$$
\end{conj}

Note that, for an Einstein metric $h$ on an $n$-manifold $N\ (n \geq 3)$,
we have that $E(h) = Y(N, [h])$.
If $E(h) > 0$, this comes from Obata's Theorem \cite{Obata},
while the nonpositive case comes from the uniqueness
(up to a normalizing constant) of Yamabe metrics.
Setting
$$Y_1 := Y(S^3) = 6\cdot\textrm{Vol}_{g_0}(S^3)^{2/3},\quad
Y_k := \frac{Y_1}{k^{2/3}}\quad \textrm{for}\ \ k \in {\mathbb N},$$ this conjecture implies, in particular,
that if $S^3/\Gamma$ is a smooth quotient by a finite group $\Gamma$ of order $k$, then $Y(S^3/\Gamma) = Y_k$.
This conjecture is {  known} to be true for any closed manifold uniformized by $T^3$ and for $S^3$.
In \cite{BN}, Bray and the second author showed that it is also true for $S^3/\Z_2=\RP^3$.

\begin{thm}[Bray--Neves]\label{main1}
Let $M$ be a closed $3$-manifold satisfying the following conditions$:$
\begin{enumerate}
\item[(i)] $M$ is not diffeomorphic to $S^3 ;$
\item[(ii)] $M$ does not contain any non-separating embedded $2$-sphere$;$
\item[(iii)] $M$ does not contain any two-sided real projective plane $\RP^2$.
\end{enumerate}
Then,
$$Y(M, [g]) \leq Y_2$$
for any conformal class $[g]$. Moreover, the equality holds if and only if
$(M, g)$ is conformally equivalent to $\RP^3$ with the metric $h_0$ of constant curvature one.
In particular, $Y(M) \leq Y_2.$
\end{thm}

\begin{rmk}
\
\begin{enumerate}
\item[(a)] The result concerning the conformal equivalence to $(\RP^3, [h_0])$
was not investigated in \cite{BN}.
However, it follows easily from the characterization of the equality condition
in the \textit{Riemannian Penrose Inequality} \cite[Section~8]{HI}.
\item[(b)] Condition (iii) in Theorem \ref{main1} is vacuous
when $M$ is orientable.
A more general condition was stated in \cite{BN}
and it includes the case where $M$ is $\RP^2\times S^1$.
The necessity of such condition is  addressed later.
\end{enumerate}
\end{rmk}

As a consequence, they obtained the following corollaries.

\begin{cor}
$$Y(\RP^3) = Y_2, \ Y(\RP^2 \times S^1) = Y_2.$$
\end{cor}
\begin{cor} The only prime $3$-manifolds with Yamabe invariant greater than $Y_2$
are $S^3$, $S^2 \times S^1$ and $S^2 \tilde{\times} S^1$
$($the nonorientable $S^2$ bundle over $S^1)$.
\end{cor}
All these manifolds have Yamabe invariant equal to $Y_1$.
The fact that $Y(S^2 \times S^1) = Y_1$ was proven independently
by Kobayashi \cite{Kobayashi} and Schoen \cite{Schoen-2}.

We now describe, in an informal way,
the idea behind the proof of Theorem \ref{main1}.
Most of the definitions can be found in Section \ref{tool} but,
for a more elaborated exposition,
the reader is encouraged to see \cite[Sections 3 and 4]{BN}.
All the statements regarding the \imc flow are proven in \cite{HI}.

The general approach was to construct a nice test function, not on $(M,g)$,
but instead on a conformally equivalent scalar-flat, \as manifold $(M-\{p\},g_{AF})$,
where $p$ is a point in $M$.
Bray and Neves observed that if this procedure is done for $(\RP^3, h_0)$,
then the optimal test function $u_0$ on the corresponding scalar-flat, \as manifold
(which is the spacial Schwarzschild manifold),
has the same level sets as the \imc flow $\{\widehat{\Sigma}_t\}_{t\geq 0}$
starting on the unique outermost minimal surface.

The main idea was to construct a test function $u$ on the \as manifold $M-\{p\}$
having the same level sets as the \imc flow $\{\Sigma_t\}_{t\geq0}$
starting on a connected component of the outermost minimal surface.
For each time $t$, $u$ assumes on $\Sigma_t$ the value
that $u_0$ assumes on $\widehat{\Sigma}_t$.

They showed next that if the Hawking quasi-local mass is nondecreasing,
which happens precisely if $\chi(\Sigma_t)\leq 2$ for all $t$, then $E(u^4g_{AF}) \leq Y_2$.
The topological assumptions on Theorem \ref{main1} assure that,
if the initial condition has Euler characteristic less or than or equal to two,
then this condition is preserved by the \imc (see Theorem \ref{geroch} and Lemma \ref{mono})
and so the result follows.

A question left open in \cite{BN} was how to compute $Y(\RP^3\#(S^2\times S^1))$ or
$Y((\RP^2\times S^1)\#(\RP^2 \times S^1))$.
The problem with the above idea is that, in these cases,
it is impossible to guarantee the monotonicity of the Hawking quasi-local mass.
The lack of monotonicity is a serious issue.
As a matter of fact, we point out that the reason why this idea does not work for $S^2 \times S^1$
(which has the Yamabe invariant $Y_1$),
is precisely because the Hawking quasi-local mass does not remain monotone
(for some interval of time, the flow might become a disconnected union of $2$-spheres).
In this paper, we find a way to overcome this difficulty.

$ $

\noindent {\bf  Acknowledgements.}
The first author would like to thank Richard Schoen and Rafe Mazzeo
for helpful discussions during his short visit to Stanford University in March 2003.
The second author is thankful to Steve Kerkchoff for many helpful discussions.
He would also like to express his gratitude to his advisor,
Richard Schoen, for all of his invaluable guidance.

\section{Main results}

We now describe the main results of this paper.
We denote the connected sum of $k$ copies of a closed manifold $M$ by $\#kM$
and the nonorientable $S^2$-bundle over $S^1$ by $S^2 \tilde{\times} S^1$.

The first result is a generalization of Theorem \ref{main1}.

\begin{thm}\label{main}
$(1)$ Let $M$ be a closed $3$-manifold satisfying the following conditions$:$
\begin{enumerate}
\item[(i)] $M$ is not diffeomorphic to $S^3 ;$
\item[(ii)] $M$ does not contain any non-separating embedded $2$-sphere$;$
\item[(iii)] $M$ does not contain
any two-sided real projective plane $\RP^2$.
\end{enumerate}
If $L$ is a closed $3$-manifold which has a non-separating $2$-sphere, then
$$Y(M\#L,[g]) < Y_2$$
for any conformal class $[g]$. In particular, $Y(M\#L) \leq Y_2.$

$(2)$ Let $M$ be a closed $3$-manifolds with a two-sided $\RP^2$,
then
$$Y(M, [g]) < Y_2$$
for any conformal class $[g]$. In particular, $Y(M) \leq Y_2$.
\end{thm}

The next corollary follows easily from Theorem \ref{main1} and Theorem~\ref{main}.
\begin{cor}
For any closed $3$-manifold $L$ which is not diffeomorphic to $S^3$, then
$$Y(\RP^3 \# L, [g]) < Y_2 \quad\mbox{and}\quad
Y((\RP^2 \times S^1) \# L, [\tilde{g}]) < Y_2$$
for any conformal classes $[g]$ and $[\tilde{g}]$, respectively.
In particular, $$ Y(\RP^3 \# L) \leq Y_2 \quad\mbox{and}\quad
Y((\RP^2 \times S^1) \# L)\leq Y_2.$$
\end{cor}

{  Given any two  manifolds} $M_1$ and $M_2$ with nonnegative Yamabe invariants,
Kobayashi~\cite{Kobayashi} proved the following fundamental inequality for the connected sum $M_1 \# M_2$
$$Y(M_1 \# M_2) \geq \min\{Y(M_1),\,Y(M_2)\}.$$
This inequality, Theorem \ref{main}, and the fact that
\begin{equation}\label{eq5}
Y(S^2 \times S^1) = Y_1 = Y(S^2 \tilde{\times} S^1)
\end{equation}
imply the next corollary.

\begin{cor}
For any nonnegative integers $k, \ell, m, n$ with $k + \ell \geq 1$,
we have
$$Y(\# k \RP^3 \# \ell (\RP^2\times S^1)
\# m (S^2 \times S^1) \# n (S^2 \tilde{\times} S^1)) = Y_2.$$
\end{cor}

Finally, we obtain the following classification.

\begin{cor}\label{improv}
The only closed $3$-manifolds with Yamabe invariant greater than $Y_2$ are
$$S^3\quad \mbox{and}\quad\#k (S^2 \times S^1) \# \ell(S^2 \tilde{\times} S^1),$$
where $k$ and $\ell$ are any nonnegative integers with $k + \ell \geq 1$.
\end{cor}
\begin{proof}
This follows from the Prime Decomposition Theorem,
which says that any closed $3$-manifold, not diffeomorphic to $S^3$,
can be decomposed as
$$M \#k (S^2 \times S^1)\#\ell (S^2 \tilde{\times}S^1) \#m Q $$
for some nonnegative integers $k, \ell$, and $m$ which is either $0$ or $1$.
Here, $M$ satisfies the conditions (ii) and (iii) of Theorem \ref{main},
and $Q$ is a connected sum of prime $3$-manifolds
each of which has a two-sided $\RP^2$ (see~\cite{Hempel}).
Theorem \ref{main} implies that $M$ is $S^3$ and that $m = 0$.
The statement follows from Kobayashi's inequality combined with \eqref{eq5}.
\end{proof}

\section{Tool box review}\label{tool}
The rest of the paper is devoted to the proof of Theorem~\ref{main},
and so we start by briefly reviewing the indispensable key results.
The main references are \cite{BN}, \cite{HI},
and \cite{LP} for the first two subsections
and \cite{Aubin-1} for the last subsection.

\subsection{Yamabe constants of \as $3$-manifolds.}\label{sec.sobo}
For any closed conformal $3$-manifold $(M, C)$ with positive Yamabe constant,
there exists a scalar-flat, \as metric $g_{AF} \in C$ on $M - \{\textrm{point}\}$.
It turns out that
this \as manifold is more convenient for studying Yamabe constants.
This has played a central role in \cite{Schoen-1}, \cite{BN},
and it will also be important in our approach.

\begin{defi}\label{defi}
A Riemannian $3$-manifold $(N, h)$ is said to be {\it asymptotically flat}
if there exists a compact set $K \subset N$ such that $N - K$ is diffeomorphic to
$\R^3 - \{|x| \leq 1\}$ \\
and that, in the coordinate chart defined by this diffeomorphism,
the metric $h = \sum_{i,j}h_{ij}(x)dx^idx^j$ satisfies the following:
$$
h_{ij} = \delta_{ij} + O(|x|^{-1}),
\quad h_{ij,k} = O(|x|^{-2}),
\quad h_{ij,kl} = O(|x|^{-3}).
$$
\end{defi}

Let $(M,g)$ be a closed $3$-manifold with positive scalar curvature.
Define the {\em conformal Laplacian} to be
$$L_{g} := - 8 \Delta_{g} + R_g,$$
and consider the normalized Green's function $G_p$ for $L_g$
with a pole at $p \in M$, that is,
$$L_g G_p = c_0\cdot \delta_p \quad \mbox{on}\quad M
\quad \mbox{and}\quad \lim_{q \rightarrow p} \mbox{dist}(p,q)G_p(q) = 1.$$
Here, $c_0 > 0$ and $\delta_p$ stand respectively for
a specific universal positive constant and the Dirac $\delta$-function at $p$.
The maximum principle implies that $G_p$ is positive on $M^*:=M-\{p\}$,
and hence $g_{AF}:=G_p^4 g$ defines a scalar-flat, \as metric on $M^*$
(see~\cite[Section~6]{LP}).

For the noncompact manifold $(M^*, g_{AF})$,
the Yamabe constant is defined by
\begin{align*}
Y(M^*, [g_{AF}]) & := \inf \left\{
\frac{\int_{M^*} 8 |\nabla u|^2 dV_{g_{AF}}}{\left(\int_{M^*} u^6 dV_{g_{AF}}\right)^{1/3} }
\:\Bigg|\: u \not\equiv 0 \mbox{ in }  C_c^{\infty}(M^*)\right\} \\
& = \inf \left\{
\frac{\int_{M^*} 8 |\nabla u|^2 dV_{g_{AF}}}{\left(\int_{M^*} u^6 dV_{g_{AF}}\right)^{1/3} }
\:\Bigg|\: u \not\equiv 0 \mbox{ in } W^{1,2}(M^*,g_{AF})\right\}.
\end{align*}
{  From} the Sobolev embedding theorem for \as manifolds, we note that $$W^{1,2}(M^*,g_{AF})\hookrightarrow
L^{6}(M^*,g_{AF}).$$ It is readily seen that $Y(M^*, [g_{AF}]) = Y(M, [g])$.

\subsection{Inverse mean curvature flow}\label{sec.imc}

Like in \cite{BN}, the inverse mean curvature flow
developed by Huisken and Ilmanen in \cite{HI}
will be a crucial tool to prove Theorem~\ref{main}.
The main difference this time is its application to
a {\it family} of \as manifolds {\it with boundary}.
We now introduce some of their terminology
and present the main results of \cite{HI} that will be necessary for our purposes.

In this subsection, $(N,\,h)$ will always denote an \as $3$-manifold with boundary $\partial N$
consisting of a union of minimal surfaces.
For a $C^1$ surface $\Sigma$ contained in $N$, $H \in L_{\l}^1(\Sigma)$
is called the {\em weak mean curvature} of $\Sigma$ provided
\begin{equation*}
\int_{\Sigma} \d_N ({X}) \,dA_h
= \int_{\Sigma} H \langle {X},\nu \rangle \,dA_h
\end{equation*}
for all vector fields $X$ with compact support on $N$,
where $\nu$ is the exterior unit normal vector.
If $\Sigma$ is smooth, $- H\nu$ coincides with the usual mean curvature of $\Sigma$.

For any surface $\Sigma$ of $(N,h)$, $\big |\Sigma\big|=\big |\Sigma\big|_h$
stands for the $2$-dimensional Hausdorff measure of $\Sigma$.
A precompact set $E$ of $N$ is called a {\it minimizing hull} {\it in $N$} if
it minimizes area on the outside, that is,
$$\big |\partial^* E\big | \leq \big|\partial^*F\big|$$
for any $F$ containing $E$ such that $F - E \subset\subset N$,
where $\partial^*F$ is the reduced boundary of $F$.

Huisken and Ilmanen \cite[Sections 3 and~6]{HI}
proved the following existence result for \imc flow.

\begin{thm}[Existence]\label{exist}
Let $E_0$ be an open precompact minimizing hull in $N$
with smooth boundary $\partial E_0$.
Then, there exists a precompact locally Lipschitz $\phi$ on $N$ satisfying$:$
\begin{enumerate}
\item[(i)] For $t \geq 0$, $\Sigma_t := \partial\{\phi < t\}$
defines an increasing family of $C^{1,\alpha}$ surfaces
$(0 < \alpha < 1)$ such that $\Sigma_0=\partial E_0;$
\item[(ii)] For almost all $t \geq 0$,
the weak mean curvature of $\Sigma_t$ is $|\nabla \phi\,|_{\Sigma_t};$
\item[(iii)] For almost all $t \geq 0$,
\begin{equation*}
|\nabla \phi\,|_{\Sigma_t} \neq 0 \quad \mbox{for almost all}\ \ x \in \Sigma_t
\end{equation*}
$($with respect to the surface measure$)$ and
\begin{equation*}
\big|\Sigma_0 \big| \, e^t \leq\big|\Sigma_t \big|\leq
(\big|\Sigma_0 \big|+\big|\partial N \big|)\, e^t \quad \mbox{for all}\ \ t \geq 0.
\end{equation*}
\end{enumerate}
\end{thm}

When $\phi$ is smooth with nonvanishing gradient, (ii) implies that
the surfaces $\{\phi = t\}_{t \geq 0}$ are a solution to \imc flow
with initial condition $\partial E_0$ because the velocity vector
for such deformation is given by
\begin{equation*}
 |\nabla \phi\,|^{-1}\nu,
\end{equation*}
where $\nu$ is the exterior unit normal to $\{\phi=t\}$.

\begin{defi}
For a compact $C^1$ surface $\Sigma$ with weak mean curvature $H$ in $L^2(\Sigma)$,
the {\it Hawking quasi-local mass} is defined to be
\begin{equation*}
m_H(\Sigma) := \sqrt{\frac{\displaystyle \big |\Sigma \big |}
{\displaystyle (16\pi)^3}}\left (16\pi-\int_{\Sigma}H^2dA_h\right).
\end{equation*}

\end{defi}
For the rest of this subsection, assume that $E_0$ is an open connected precompact {minimizing hull in $N$ with}
smooth boundary $\partial E_0$, and that $\{\Sigma_t\}_{t \geq 0}$ is a weak solution for the \imc flow with
initial condition $\Sigma_0 = \partial E_0$.

Under the condition that $(N,\,h)$ has nonnegative scalar curvature,
we have the following remarkable property for the  weak \imc flow \cite{HI}.

\begin{thm}[Geroch Monotonicity]\label{geroch}
Let $(N,h)$ be an \as $3$-manifold with nonnegative scalar curvature.
The Hawking quasi-local mass $m_H(\Sigma_t)$ is a nondecreasing function of $t$
provided $\chi (\Sigma_t) \le 2$ for all $t$.
\end{thm}

The monotonicity of the Hawking quasi-local mass will be of major importance
in proving Theorem~\ref{main}.
The next Lemma gives topological conditions on $N$ that guarantee this monotonicity.
\begin{lemm}\label{mono}
Assume that besides the assumption in Theorem~\ref{geroch},
$N$ does not contain any non-separating embedded $2$-sphere
and any two-sided real projective plane.
Then,
$$\chi (\Sigma_t) \leq 2$$
for all $t \geq 0$, provided that $\partial E_0$ is connected.
\end{lemm}
\begin{proof}
Let $\phi$ be the function defining the weak solution $\{\Sigma_t\}_{t \geq 0}$. According to \cite[Sections~4
and~6]{HI}, { the connectedness of $E_0$} implies that both $\{\phi < t\}$ and $\{\phi
> t\}$ remain connected for all $t$.

It is also proven in  \cite{HI} that, for  each $t >0$,
$\Sigma_t$ can be approximated by earlier surfaces $\Sigma_s$
for which $\nabla \phi$ never vanishes.
Hence, without loss of generality, we can assume that
\begin{equation*}
\Sigma_t = \{\phi = t\} = \partial \{\phi > t\} = \partial \{\phi < t\}.
\end{equation*}
The connectedness of the two sets on the right of the above expression
implies that $\Sigma_t$ is two-sided.
We now show that $\chi (\Sigma_t) > 2$ is impossible.

If $\Sigma_t$ contains a $2$-sphere,
then this $2$-sphere has to bound a connected region,
and hence $\chi (\Sigma_t) = 2$ because $\{\phi < t\}$ is connected.
If $\Sigma_t$ contains a real projective plane,
then it has to be one-sided.
This contradicts that $\Sigma_t$ is two-sided.
\end{proof}

\subsection{Aubin's lemma}
The next Lemma, due to Aubin \cite{Aubin-1},
is a key tool to prove Theorem~\ref{main}.

\begin{lemm}[Aubin]\label{aubin}
Let $(M, g)$ be a closed $n$-manifold with positive Yamabe constant
and $(M_k, g_k)$ a non-trivial finite $k$-fold Riemannian covering $(i.e., k \geq 2)$.
Then,
$$Y(M, [g]) < Y(M_k, [g_k]).$$
\end{lemm}

\begin{rmk}
{Aubin's proof assumes that $M_k$ is a normal covering of $M$ (cf.~\cite[Theorem~6]{Aubin-1}). We have to adapt
its proof in order to remove this condition.}
\end{rmk}

\begin{proof}
Let $u$ {be} a normalized positive function that minimizes the Yamabe quotient $Q_{g_k}$, that is,
$$Q_{g_k}(u) = E(u^{4/(n-2)}g_k) = Y(M_k, [g_k])\quad
\mbox{and}\quad \int_{M_k}u^{2n/(n-2)} \,dV_{g_k} = 1.$$ For any $x \in M_k$, {set
$$\{x_1, \cdots, x_k\}:= \mathcal{P}^{-1}\big{(}\mathcal{P}(x)\big{)}\quad\mbox{with}\quad x_1 = x,$$
where $\mathcal{P}$ denotes the $k$-fold covering map.} Define positive functions $v$ and $v^{<p>}$ (for each $p
> 0$) on $M_k$ by
$$ v(x) := \sum_{i=1}^k u(x_i)\quad\mbox{and}\quad v^{<p>}(x) := \sum_{i=1}^k u(x_i)^p.$$
For any {\it evenly covered} open set $U\ (\subset M)$ for $\mathcal{P}$,{ set
$$\{U_1, \cdots, U_k\} := \mathcal{P}^{-1}(U)\ (\subset M_k).$$}
{ Moreover, all the sets are isometric and so we can find, for instance, $k$ isometries
$$\gamma_i : U_1 \longrightarrow U_i\quad i = 1, \cdots, k$$} satisfying $\mathcal{P} \circ \gamma_i =
\mathcal{P}$. {  As a result, it is possible to} express $v$ and $v^{<p>}$ on $U_1$ as
$$v = \sum_{i=1}^k u \circ \gamma_i,\quad
v^{<p>} := \sum_{i=1}^k (u \circ \gamma_i)^p,$$ and { hence} $v, v^{<p>} \in C^{\infty}(M_k)$.

{We can use the local expression described above in order find  $k$ distinct isometries
$$\gamma_i : M_k - \mathcal{S} \longrightarrow M_k - \mathcal{S},\quad i = 1, \cdots, k,$$
where $\mathcal{S}$ is a piecewise smooth compact ($n -1$)-submanifold  such that
$\mathcal{P}^{-1}\left(\mathcal{P}(\mathcal{S})\right)=\mathcal{S}$.
 Hence, if we set $u_i:=u\circ\gamma_i$, we have that on $M_k - \mathcal{S}$


$$v = \sum_{i=1}^k u_i \quad\mbox{and}\quad v^{<p>} = \sum_{i=1}^k u_i^p,
$$ }
from which it follows that
\begin{equation*}
\int_{M_k}v^{<p>}\,dV_{g_k} = k\ \int_{M_k}u^p \,dV_{g_k}
\end{equation*}
and $$L_k v = Y(M_k, [g_k])v^{<(n+2)/(n-2)>} = Y(M_k, [g_k]) \sum_{i=1}^k u_i^{(n+2)/(n-2)},$$
 where $L_k := - \frac{4(n-1)}{n-2} \Delta_{g_k} + R_{g_k}$ is
the conformal Laplacian with respect to $g_k$.

Define $v_0$ to be the function on $M$ whose lift to $M_k$ is $v$.
Then,
\begin{equation}\label{aub1}
\begin{split}
Q_{g}(v_0) & = k^{-2/n}Q_{g_k}(v)\\
& = k^{-2/n} \frac{\displaystyle \int_{M_k}v (L_k v) \,dV_{g_k}}
{\left( \displaystyle \int_{M_k} v^{2n/(n-2)} \,dV_{g_k} \right)^{(n-2)/n}}\\
& = k^{-2/n} Y(M_k, [g_k]) \frac{\displaystyle \int_{M_k} v \left(\sum_{i=1}^k u_i^{(n+2)/(n-2)}\right)
\,dV_{g_k}} {\left( \displaystyle \int_{M_k} v^{2n/(n-2)} \,dV_{g_k} \right)^{(n-2)/n}}.
\end{split}
\end{equation}
Using H\"older's inequality twice and the following strict inequality
$$\left(\sum_{i=1}^k a_i^p\right)^{1/p} < \sum_{i=1}^k a_i\qquad
\textrm{for}\quad k \geq 2,\ p> 1,\
a_i > 0\ (i = 1, ..., k), $$
we obtain that

\begin{equation}\label{aub2}
\begin{split}
& \displaystyle \int_{M_k} v  \left(\sum_{i=1}^k u_i^{(n+2)/(n-2)}\right) dV_{g_k}\\
& \leq  \displaystyle \int_{M_k} v \left( \sum_{i=1}^k u_i^{2n/(n-2)} \right)^{2/n}
\left( \sum_{i=1}^k u_i^{n/(n-2)} \right)^{(n-2)/n}dV_{g_k}\\
& <  \displaystyle \int_{M_k} v^2
\left( \sum_{i=1}^k u_i^{2n/(n-2)} \right)^{2/n}dV_{g_k} \\
& \leq \left(\displaystyle\int_{M_k} v^{2n/(n-2)} \,dV_{g_k} \right)^{(n-2)/n}
\left(\displaystyle\int_{M_k}\sum_{i=1}^k u_i^{2n/(n-2)} dV_{g_k} \right)^{2/n}\\
& = \left(\displaystyle\int_{M_k} v^{2n/(n-2)} \,dV_{g_k} \right)^{(n-2)/n}
\left(\displaystyle\int_{M_k} v^{\langle 2n/(n-2)\rangle} dV_{g_k} \right)^{2/n}\\
& =  k^{2/n} \left(\displaystyle\int_{M_k} v^{2n/(n-2)} \,dV_{g_k} \right)^{(n-2)/n}.
\end{split}
\end{equation}
It  {  follows} from \eqref{aub1} and \eqref{aub2} that
$$Q_g(v_0) < Y(M_k, [g_k])$$
and hence
$$Y(M, [g])<Y(M_k, [g_k]).$$
\end{proof}

\section{Overall strategy}\label{heuristic}
In this section we describe our approach to Theorem~\ref{main}-$(1)$
in the simple case that $L = S^2\times S^1$.
Let
$$N :=M \# (S^2 \times S^1),$$
where $M$ is a closed $3$-manifold satisfying conditions (i), (ii), and (iii) in Theorem~\ref{main}.
Without loss of generality, we may assume that $Y(N, [g]) > 0$
and that $g$ is a Yamabe metric with Vol$_g(N) = 1$.
Note that $R_g \equiv Y(N, [g]) >0$.

Consider Riemannian coverings of $(N, g)$ defined in the following way.
Take a smooth loop $c$ in $S^2 \times S^1$ whose homotopy class $\langle c\rangle$ generates $\pi_1(S^2\times S^1)$.
We regard the loop $c$ as a loop in $N$.
Let $(N_k, g_k)$ and $(N_{\infty}, g_{\infty})$ denote, respectively,
the $k$-fold Riemannian covering and the normal infinite Riemannian covering associated to $\langle c \rangle$.
Note that, topologically,
$$N_k = \# kM \# (S^2 \times S^1), \qquad
N_{\infty} = \#_1^{\infty} M \# (S^2 \times \R),$$ and that the covering $N_{\infty} \rightarrow N_k$ is also
normal. From Aubin's Lemma~\ref{aubin} we have that, for all $k\geq 2$,
$$Y(N, [g]) < Y(N_k, [g_k]),$$
and this implies that
$$Y(N, [g]) < \limsup_{k \to \infty} Y(N_k, [g_k]).$$
The strict inequality comes from the fact that
we can choose an increasing subsequence of $\{Y(N_k, [g_k])\}_{k\geq1}$,
for instance, $\{Y(N_{2^m}, [g_{2^m}])\}_{m\geq1}$.
Moreover, modifying the techniques in \cite[Section~6]{Kobayashi}
and \cite[Sections~5 and~6]{AB},
we can prove the following
$$
\lim_{k \to \infty} Y(N_k, [g_k]) = Y(N_{\infty}, [g_{\infty}]).
$$
However, the above strict inequality is sufficient
for our proof that $Y(N, [g]) < Y_2$.

Unlike each $N_k$, the noncompact manifold $N_{\infty}$ satisfies the conditions (i), (ii), and (iii) of
Theorem \ref{main1}, and hence, at least heuristically, this suggests that
$$\limsup_{k \to \infty} Y(N_k, [g_k]) = Y(N_{\infty}, [g_{\infty}]) \leq Y_2.$$
It should be also pointed out that the scalar-flat, \as manifold conformally equivalent to $(N_{\infty} -
\{\mbox{point}\}, g_{\infty})$, will have two singularities created by the two ends of $N_{\infty}$. Therefore,
we can not apply the \imc flow technique directly to this manifold. Instead, we apply it to a family of \as
manifolds without singularities { (but with boundary), arising from the manifolds $(N_k^*, g_{k,AF})$ that we
describe next.}

Fix a point $p_{\infty} \in N_{\infty}$
and call its projection to each $N_k$ also by $p_{\infty}$.
Let $G_k$ be the normalized Green's function on $N_k$
for the conformal Laplacian $L_{g_k}$ with the pole at $p_{\infty}$.
Consider the following scalar-flat, \as metrics
$$g_{k,AF} := G_k^4g_k \quad \mbox{on}\quad N_k^* := N_k - \{p_{\infty}\}.$$

{ Using inverse mean curvature flow techniques, we will show}
\begin{thm}\label{main3}
$$ \limsup_{k \to \infty} Y(N_k^*, [g_k]) \leq Y_2. $$
\end{thm}

The proof of Theorem \ref{main3} is given in the next two sections.
Some extra difficulties arise for the general case
because the infinite Riemannian covering will have infinitely many ends.
Fortunately, the method developed here and in Sections~5, 6 is robust enough to handle those cases.
The proof will be given in Section~\ref{general}.

\section{Proof of main theorem: simple case}\label{proof}

Under the same assumptions of Section~\ref{heuristic},
consider the normalized minimal positive Green's function $G_{\infty}$
on $(N_{\infty}, g_{\infty})$ for the conformal Laplacian $L_{g_{\infty}}$
with the pole at $p_{\infty}$ (cf.~\cite[Section~1]{SY-3}, \cite[Section~6]{AB}).
This function is defined on $N^*_{\infty}:=N_{\infty}-\{p_{\infty}\}$
and it is essential to control the geometry of $(N_k^*, g_{k,AF})$.
All the statements in this section regarding these manifolds will follow
from the good understanding we have for $G_{\infty}$.

We start by collecting the properties $G_{\infty}$ and $G_k$ in Lemma~\ref{greens}.
After this, we describe the main idea behind the proof of Theorem~\ref{main3} and state four auxiliary lemmas,
whose proofs will be given in the next section.
We then prove Theorem~\ref{main3}.

For reasons of simplicity, we set the following notation for the rest of the paper. Any quantity $Q$ that
depends on $g_{k,AF}$ or $g_k$ will be denoted by $Q_{k,AF}$ or $Q_k${ respectively}. If it is obvious from the
context which metric is being considered, we will just denote it by $Q$.

Let $\gamma$ be the deck transformation of $N_{\infty}$
corresponding to the homotopy class of $\langle c \rangle \in \pi_1(N)$,
that is,
$$N_k = N_{\infty}/\{\langle \gamma^{jk}\rangle\,|\,j\in\Z \}\quad\mbox{for}\quad k \geq 1.$$
For simplicity, we often identify a set in $N_{\infty}$ with its projection to each $N_k$ and a function on
$N_k$ with its lift to $N_{\infty}$. For instance, $G_k$ denotes both the normalized Green's function on $N_k$
and its lift to $N_{\infty}$.

{ Standard arguments developed in \cite{SY-3}, \cite{Haber}, and \cite{AB}, imply the following lemma regarding
the functions $G_k$ and $G_{\infty}$.} Here, note that each $g_k$ and $g_{\infty}$ are metrics of positive
constant scalar curvature $Y(N, [g])
> 0$.

\begin{lemm}\label{greens}
\mbox{}
\begin{enumerate}
\item[(i)] For all $k\geq 1$,
$$G_k = \sum_{j \in \Z} G_{\infty}\circ\gamma^{jk}\quad\mbox{on }N_{\infty};$$
\item[(ii)] For all $k\geq 2$,
$$0 < G_{\infty} < G_k < G_1,$$ and
$G_k$ converges uniformly in $C^{\ell}$ $($ for every $\ell \geq 1)$ to
$G_{\infty}$ on every compact subset of $N_{\infty}^*;$
\item[(iii)] For every compact set $K$ of $N_{\infty}^*$,
$$\lim_{j \to \infty} \sup \big \{G_{\infty}(\gamma^j (x))\:|\: x \in K \big \} = 0;$$
\item[(iv)] For any open set $O$ containing $p_{\infty}$,
there exists a constant $L$ independent of $k$ such that, for all $k\geq 1$,
$$|\nabla G_{k}|_k\le L G_{k}
      \quad\mbox{on }N_{k}-O;$$
\item[(v)] $$\int_{N_{\infty}}G_{\infty}\,dV_{\infty} < \infty. $$
\end{enumerate}
\end{lemm}

Take a smooth connected embedded $2$-sphere $S$ in $N_{\infty}$ with $p_{\infty} \in S$
such that the projection from $S$ into $N$ is injective
and its image in $N$ intersects transversely $c$ at one point.
Then, consider the embedded $2$-spheres
$$F_k := \gamma^{[k/2]}(S) \sqcup \gamma^{[k/2] - k}(S)$$
in $N_{\infty}$ for each $k \geq 2$.
Denote the projection of $F_k\ (\subset N_{\infty})$
to each $N_k$ also by $F_k$.
Note that $F_k$ is connected and that
$$\lim_{k \to \infty} \mbox{dist}_{g_k}(p_{\infty}, F_k) = \infty.$$

The topological condition on $M$ implies that
all non-separating $2$-spheres of $N_k^*$ are cobordant to each other.
Such $2$-spheres correspond to a `` neck'' on $N^*_k$.
The properties (ii) and (iii) say that $(N_k^*, g_{k,AF})$ has,
for $k$ large enough, a very small ``neck''.
This is the content of the following Lemma.

\begin{lemm}\label{small}
For every $k$, there exists an area-minimizing, non-separating,
connected embedded $2$-sphere
$S_k$ of $(N_k^*, g_{k,AF})$ such that
$$\lim_{k \to \infty}\big |S_k \big|_{k,AF} = 0,$$
where $\big |S_k \big|_{k,AF}$ denotes the area of $S_k$ with respect to $g_{k,AF}$.
\end{lemm}

Set $X_k$ to be the metric completion of
$$N_k^* - S_k = N_k - (\{p_{\infty}\} \cup S_k),$$
and denote the smooth metric extension of $g_{k,AF}$ to $X_k$ by the same symbol.
Then, $(X_k,g_{k,AF})$ becomes an \as manifold with minimal boundary consisting of two $2$-spheres.
This corresponds to ``cutting out the small neck'' of $N_k^*$.
The advantage of working with $X_k$ rather than $N_k^*$ is that, as can be easily seen,
any embedded $2$-sphere of $X_k$ bounds a connected region.
Therefore, by Theorem~\ref{geroch} and Lemma~\ref{mono}, the monotonicity of the Hawking quasi-local mass
is preserved for each solution of \imc flow on $(X_k, g_{k,AF})$,
provided that its initial condition is the connected boundary of a connected region.
Taking into account that the area $|\partial X_k|_{k,AF}$ of $\partial X_k$ goes to zero as $k$ goes to infinity,
we expect to construct a nice family of test functions $u_k$ with
$$\limsup_{k\to\infty}Q_{g_{k,AF}}(u_k)\leq Y_2.$$

Using the \imc flow techniques in the same spirit as in \cite{BN},
we have the following result.

\begin{lemm}\label{prop}
Let $(X, h)$ be an \as $3$-manifold with minimal boundary $\partial X$.
Assume further there exists a weak solution for the \imc flow,
with initial condition a minimal surface $\Sigma_0$,
that keeps the Hawking quasi-local mass nondecreasing.
Then, there exists a function $v$ in $W^{1,2}(X,h)$ such that
$0 \leq v \leq 1$ on $X$ and
$$\frac{\int_X 8 |\nabla v|^2 dV_h}{\left(\int_X v^6 dV_h\right)^{1/3} }
\le Y_2 + \theta \left(\big|\partial X \big|\big |\Sigma_0 \big|^{-1}\right),$$
where $\theta$ is a continuous function independent of $(X, h)$
and vanishing at zero.
\end{lemm}

An {\em outermost minimal surface} of an \as $3$-manifold $(X,h)$
is a minimal surface (possibly disconnected)
which encloses all other compact minimal surfaces.
Since each $N^*_k$ is not diffeomorphic to $\R^3$,
a result due to Meeks, Simon and Yau \cite{MSY} implies that
$(X_k,g_{k,AF})$ has a unique outermost minimal surface $\Omega$
which is the disjoint union of a finite number of weakly embedded $2$-spheres
(i.e., limits of uniformly smooth embedded $2$-spheres).
Moreover, the region exterior to the outermost minimal surface
is $\R^3$ minus a finite number of closed balls (see \cite[Section~4]{HI} for details).
The next Lemma guarantees that
one of the connected components of the outermost minimal surface
has area bounded below independently of $k$.

\begin{lemm}\label{essential}
The outermost minimal surface of each $(X_k, g_{k,AF})$ has a connected component $\Sigma(k)$
with area bounded below independently of $k$.
\end{lemm}

\begin{rmk}\label{tech}
The surface $\Sigma(k)$ is a minimizing hull
because it is contained in the outermost minimal surface.
If $\Sigma(k)$ is a $2$-sphere,
then it is the boundary of a precompact connected open set,
and hence it can be used as an initial condition for the \imc flow (see Theorem~\ref{exist}).
If $\Sigma(k)$ is a one-sided $\RP^2$,
then it is the limit of uniformly smooth embedded $2$-spheres
that are minimizing hulls bounding a precompact open connected set.
We then take the limit of the weak solutions to \imc flow starting at each of these $2$-spheres,
in order to get a weak solution to the \imc flow with initial condition $\Sigma(k)$.
This can be done rigorously using the compactness theorem
and the a priori bounds proved in \cite[Sections~2 and 3]{HI}.
\end{rmk}

It can be easily seen that Theorem~\ref{main3} is false
when $M$ is either $S^3$ or $S^2\times S^1$,
in which case $N$ is $S^2 \times S^1$ or
$(S^2 \times S^1) \# (S^2 \times S^1)$, respectively.
The technical reason why this proof doesn't work for those cases
is because Lemma \ref{essential} doesn't hold.
It is worthwhile to point out that, among the auxiliary lemmas we state in this section,
Lemma~\ref{essential} is the only lemma
where the topological restrictions on $M$ play a role.

The function $v$ given by Lemma \ref{prop} may not vanish on the boundary.
Therefore, if we want to apply this result
to our sequence of \as manifolds $(X_k,\,g_{k,AF})$,
we need to multiply the resulting functions
by some cut-off function vanishing on the boundary $\partial X_k$.
Otherwise, we could not extend it to be a valid test function on $N_k^*$.
This is a delicate issue because, in general,
the multiplication by a cut-off function makes the Yamabe quotient increase.
Fortunately, using Lemma~\ref{greens} (iv) and (v),
we can prove the following result.

\begin{lemm}\label{greens_small}
For every $\varepsilon > 0$,
there exists a sequence of Lipschitz functions
$\eta_k$ $(0 \leq \eta_k \leq 1)$ such that
$$\eta_k \equiv 0\quad\mbox{on}\ \partial X_k,
\qquad \eta_k\equiv 1\quad \mbox{outside a compact set of}\ X_k,$$
and for all $k$ sufficiently large
$$\int_{X_k}|\nabla \eta_k|^2\,dV_{k,AF} \leq \varepsilon, \quad
\vol_{k,AF}\Big(\{\,x \in X_k\:|\:\nabla \eta_k\neq 0\,\}\Big) \leq \varepsilon.$$
\end{lemm}

We can now prove Theorem~\ref{main3}.
\begin{proof}[\bf {Proof of Theorem \ref{main3}}]

The proof follows if we find a sequence of test functions $u_k$
vanishing on $X_k$ such that
$$\limsup_{k \to \infty}
\frac{\displaystyle \int_{X_k} 8|\nabla u_k|^2dV_{k,AF}}
{\displaystyle \left(\int_{X_k} u_k^6 dV_{k,AF}\right)^{1/3} } \leq Y_2 $$

Consider the  minimal surface $\Sigma(k)$ given by Lemma~\ref{essential}. There exists a solution to the \imc
flow with initial condition $\Sigma(k)$ (Remark~\ref{tech}) and the topological condition on $X_k$ implies the
monotonicity of the Hawking quasi-local mass along the solution. From Lemma~\ref{prop}, we have the existence,
for all $k$, of $v_k$ $(0\leq v_k \leq 1)$ in $W^{1,2}(X_k, g_{k,AF})$ such that
$$\frac{\displaystyle \int_{X_k} 8 |\nabla v_k|^2dV_{k,AF}}
{\displaystyle \left(\int_{X_k} v_k^6 dV_{k,AF}\right)^{1/3} }
\leq Y_2 + \theta \left(\big|\partial X_k \big|_{k,AF}
\big|\Sigma(k) \big|^{-1}_{k,AF}\right).$$

Take any small $\varepsilon > 0$
and let $\eta_k\ (0 \leq \eta_k \leq 1)$ denote the Lipschitz function
given in Lemma~\ref{greens_small} for all sufficiently large $k$.
Consider the admissible nonnegative function in $W^{1,2}(N_k^*,\,g_{k,AF})$ given by
$$u_k := \eta_k v_k.$$
We then have

\begin{multline*}
\int_{X_k}8|\nabla u_k|^2 \,dV_{k,AF}
\leq  8\ \Bigg{\{}\int_{X_k}|\nabla v_k|^2 \,dV_{k,AF} \\
+2\sqrt{ \int_{X_k}|\nabla v_k|^2\, dV_{k,AF}}\sqrt{ \int_{X_k}|\nabla \eta_k|^2 \,dV_{k,AF}}+
\int_{X_k}|\nabla \eta_k|^2\, dV_{k,AF}\Bigg{\}}\\
 \leq q_k \left ( \int_{X_k}v_k^6\,dV_{k,AF}\right)^{1/3} +
4\sqrt{2\varepsilon q_k}\left (\int_{X_k}v_k^6\,dV_{k,AF}\right)^{1/6}+8\varepsilon,
\end{multline*}
where $$q_k:= Y_2 + \theta \Big(\big|\partial X_k \big|_{k,AF}
\big|\Sigma(k) \big|^{-1}_{k,AF}\Big)>0.$$

Following the proof of Lemma \ref{prop} we see that
$$\int_{X_k}v_k^6\,dV_{k,AF} \geq
\big|\Sigma(k)\big|_{k,AF}^{3/2}\,
\Big(C_0 +
\theta \big(\big|\partial X_k \big|_{k,AF}\big|\Sigma(k) \big|^{-1}_{k,AF}\big)\Big)
\geq C>0 $$
for some  positive constants $C_0, C$ independent of $k$,
where the second inequality comes from combining Lemma~\ref{small}
with the lower bound on $\big|\Sigma(k)\big|_{k,AF}$ (Lemma~\ref{essential}),
and with the continuity of $\theta$.
Hence, using the properties of $\eta_k$ given by Lemma~\ref{greens_small},
we obtain
\begin{align*}
\frac{\displaystyle \int_{X_k} v_k^6 \,dV_{k,AF}}{\displaystyle
  \int_{X_k} u_k^6\,dV_{k,AF}}
 &\leq  \frac{\displaystyle \int_{X_k}
  v_k^6 \,dV_{k,AF}}{\displaystyle
  \int_{X_k}v_k^6\,dV_{k,AF}-\vol_{k,AF}
  \Big(\{\,x \in X_k\:|\:\nabla \eta_k\neq 0\,\}\Big)}\\
& \leq  \frac{C}{C-\varepsilon}
\end{align*}
for some positive constant $C$ independent of $k$, and therefore
$$ \frac{\displaystyle \int_{X_k}8|\nabla u_k |^2\, dV_{k,AF} }
{\displaystyle \left( \int_{X_k} u_k^6\, dV_{k,AF} \right)^{1/3}} \leq q_k
\left(\frac{C}{C-\varepsilon}\right)^{1/3} +
4\sqrt{2q_k}\frac{\sqrt{\varepsilon}\,C^{1/6}}{(C-\varepsilon)^{1/3}} +
\frac{8\varepsilon}{(C-\varepsilon)^{1/3}}.$$ Combining Lemma~\ref{essential} with the fact that $\big|\partial
X_k \big|_{k,AF}$ converges to zero (Lemma~\ref{small}), we have
$$\lim_{k \to \infty} q_k = Y_2.$$
The proof follows  by letting $\varepsilon$ go to zero.
\end{proof}

\section{Proofs of auxiliary lemmas}\label{proof2}

We first prove Lemma~\ref{small}, then Lemma~\ref{prop},
Lemma~\ref{essential}, and finally Lemma~\ref{greens_small}.
\begin{proof}[{\bf Proof of Lemma \ref{small}}]
Consider the same smooth connected embedded $2$-spheres $S$ in $N_{\infty}$ and $F_k$ in $N_k^*$ that were
constructed right before {  stating} Lemma~\ref{small}. Then
$$\big | \gamma^j(S) \big |_{{\infty,AF}}
 = \int_S G^4_{\infty}\circ \gamma^j\,dA_{\infty}. $$
Hence, from Lemma~\ref{greens} (ii) and (iii),
$$\lim_{k \to \infty} \big| F_k \big |_{{k,AF}} = 0. $$
{ Next}, we minimize area in the isotopy class of $F_k$. The theory developed in \cite{MSY} says that the area
minimizer exists and is a smooth disjoint union of $2$-spheres and one-sided real projective planes. Moreover,
according to \cite[Section~3]{MSY}, the connected components of the area minimizer can be joined by arbitrarily
thin tubes so that the resulting surface is isotopic to $F_k$. Hence, the area minimizer has to contain a smooth
minimal $2$-sphere $S_k$ that is non-separating.

The desired assertion follows from the fact that
$\big|S_k\big|_{k,AF} \leq \big|F_k\big|_{k,AF}$.

\end{proof}

\begin{proof}[{\bf Proof of Lemma \ref{prop}}]
We use the same idea of \cite{BN}
and construct a test function $v$
whose level sets coincide with the weak solution of the \imc flow.
This time we have to take into account that $X$ has a minimal boundary $\partial X$.
Let $\phi$ be the weak solution to the \imc flow with initial condition
$\Sigma_0$ and set $\Sigma_t := \partial\{\phi < t\}$.
Like in \cite{BN}, $v$ is defined to be
\begin{equation*}
v(x):=\left\{
  \begin{array}{r@{\quad\mbox{if  }}l}
   f(0) &  \phi(x) \leq 0, \\
   f(\phi(x)) & \phi(x) > 0,
  \end{array}
  \right.
\end{equation*}
where $$f(t) = \frac{1}{\sqrt{2e^t - e^{t/2}}}.$$
After a simple computation using the coarea formula \cite[Section~6]{BN},
we get
\begin{align}
\int_{X}|\nabla v|^2 dV_h
& =  \int_0^{\infty} {f'(t)}^2 \int_{\Sigma_t} H dA_h\, dt, \label{f1}\\
\int_{X} v ^6\, dV_h
& \geq  \int_0^{\infty}f(t) ^ 6 \,
\big |\Sigma_t \big|^2  \left (\int_{\Sigma_t } H \,dA_h \right )^{-1} dt.\label{f2}
\end{align}
Next, we find an upper bound for $\int_{\Sigma_t } H\, dA_h$.
\begin{lemm}\label{aqui}
For all $t \geq 0$,
\begin{equation*}
\int_{\Sigma_t } H dA_h  \leq
\sqrt{16\pi\big|\Sigma_0\big|(e^t - e^{t/2})}+\sqrt{16\pi\big|\partial X \big |e^t}.
\end{equation*}
\end{lemm}
\begin{proof}
Using the monotonicity of the Hawking quasi-local mass, we have
\begin{align*}
\sqrt{\frac{ \displaystyle \big |\Sigma_t \big | }{\displaystyle (16\pi)^3} }
\left(16\pi-\int_{\Sigma_t}H^2 \, dA_h\right)
&= m_H(\Sigma_t) \\
& \geq  m_H(\Sigma_0) \\
&= \sqrt{\frac{\displaystyle \big |\Sigma_0 \big|}{ \displaystyle 16\pi}}.
\end{align*}
Hence, it follows from H\"{o}lder's inequality that
\begin{align*}
\int_{\Sigma_t } H dA_h & \leq  \sqrt{\big|\Sigma_t \big|\int_{\Sigma_t}H^2 \, dA_h}, \\
& \leq
\sqrt{16\pi\left(\big| \Sigma_t \big|-\sqrt{\big|\Sigma_t\big|\big|\Sigma_0\big|}\right)},
\end{align*}
and thus, by Theorem~\ref{exist} (iii),
\begin{align*}
\int_{\Sigma_t } H dA_h  &\leq
\sqrt{16 \pi \big|\Sigma_0 \big|(e^t - e^{t/2})+16 \pi \big |\partial X\big| e^t}\\
&\leq  \sqrt{16\pi\big|\Sigma_0\big|(e^t - e^{t/2})}+\sqrt{16\pi\big|\partial X\big|e^t}.
\end{align*}
\end{proof}

Combining Lemma \ref{aqui} first with \eqref{f1} and then with \eqref{f2},
we get
$$\int_{X}|\nabla v|^2 dV_h \leq
\sqrt{16\pi\big|\Sigma_0\big|}\left( \int_0^{\infty} f'(t)^2 \sqrt{e^{t}-e^{t/2}}\,dt
+ C\big|\partial X \big|\big|\Sigma_0 \big|^{-1}\right),$$
for some positive constant $C$.
Then, Theorem~\ref{exist} (iii) implies the following
\begin{align*}
\int_{X} v ^6\, dV_h &\geq  \int_0^\infty \frac{f(t)^6 \big|\Sigma_t\big|^2}{\sqrt{16\pi\big|\Sigma_0\big|(e^t - e^{t/2})}+\sqrt{16\pi\big|\partial X\big|e^t}}\,dt \\
&\geq  \big|\Sigma_0\big|^{3/2}\int_0^\infty \frac{f(t)^6 e^{2t}}{\sqrt{16\pi(e^t - e^{t/2})} + \sqrt{16\pi\big|\partial X\big|\big|\Sigma_0\big|^{-1}e^t}}dt\\
& =  \big|\Sigma_0\big|^{3/2}\left(\int_0^\infty \frac{f(t)^6 e^{2t}}{\sqrt{16\pi(e^t - e^{t/2})}}dt + {\theta}\left(\big|\partial X \big|\big|\Sigma_0 \big|^{-1}\right)\right),
\end{align*}
where $\theta$ is some function independent of the metric.
It can be easily seen that $\theta$ is continuous and vanishes
when $\big|\partial X\big|\big|\Sigma_0\big|^{-1}$ is zero.
According to \cite{BN}, we know that
$$\frac{\displaystyle (16\pi)^{2/3}\int_0^{\infty} f'(t)^2 \sqrt{e^{t}-e^{t/2}}\,dt}{\displaystyle  \left ( \int_0^\infty
f(t)^6 e^{2t}\,(e^{t}-e^{t/2})^{-1/2}dt \right )^{1/3}} =
\frac{\displaystyle Y_2}{\displaystyle 8},$$
and hence
\begin{align*}
\frac{\displaystyle \int_{X}|\nabla v |^2 dV_h }{\displaystyle \left( \int_{X} v^6 dV_h \right)  ^{1/3}}  &\leq  \frac{\displaystyle (16\pi)^{2/3}\int_0^{\infty} f'(t)^2 \sqrt{e^{t}-e^{t/2}}\,dt}{\displaystyle  \left ( \int_0^\infty
f(t)^6 e^{2t}\,(e^{t}-e^{t/2})^{-1/2}dt \right )^{1/3}}+ {\theta}\left(\big|\partial X \big|\big|\Sigma_0 \big|^{-1}\right)\\
&= \frac{Y_2}{8}+{\theta}\left(\big|\partial X \big|\big|\Sigma_0 \big|^{-1}\right),
\end{align*}
where ${\theta}$ is a continuous function vanishing at zero and independent of the metric $h$.
 \end{proof}

\begin{proof}[{\bf Proof of Lemma \ref{essential}}]

We argue by contradiction.
Suppose that there exists an infinite number of
\as $3$-manifolds $(X_k,\,g_{k,AF})$ (which we still index by $k$),
such that the area of every connected component of the outermost minimal surface
$\Omega_k\ (\subset X_k)$ goes to zero when $k$ goes to infinity.

We show the existence of a connected component $\Omega_{k,0}$ of $\Omega_k$
for which there exist two subsets $A_k,\ B_k$ of  $X_k$,
both intersecting $\Omega_{k,0}$,
and such that the distance between $A_k$ and $B_k$ (with respect to $g_{k,AF}$)
is bounded from below independently of $k$.
Thus,  $\Omega_{k,0}$ will have, for $k$ large enough,
very small area and diameter bounded from below.
Moreover, it will be apparent from the construction of $A_k$ and $B_k$
that the complement of both sets in $X_k$
has bounded sectional curvature independently of $k$.
A simple application of the monotonicity formula shows that
$\Omega_{k,0}$ cannot have very small area
(cf.~\cite{Schoen-0-1}, \cite{Schoen-0-2}).

First, we construct $A_k$.
Denote by $\widehat{N}\ (\subset N_{\infty})$
the closure of a fundamental domain of $N$ containing $p_{\infty}$ and take
a positive constant $\alpha$ such that
$$ A_k := G_{\infty}^{-1}[\alpha, + \infty)$$
contains $\widehat{N}-\{p_{\infty}\}$. Note that this set cannot be diffeomorphic to $\R^3$ minus a finite number of balls.
For all $k$ sufficiently large,
this set projects injectively on $X_k$ and so
we can identify $A_k$ with its projection on $X_k$.
Since the region exterior to $\Omega_k$ is topologically trivial,
we get that one of its connected components ${\Omega}_{k,0}$ must intersect $A_k$.

Define $$B_k := G_{k}^{-1}[0, \alpha/2].$$
We now show that, for all $k$  sufficiently large,
$B_k$ intersects $\Omega_{k,0}$.
Consider
$$s_k := \inf\left\{ \big|\Sigma\big|_{k}\: \Big|\:\ \Sigma
\mbox{ is a homotopically non-trivial surface of } N_{k} \right\},$$
and note that $0 <s_1 \leq s_k$ for all $k$.
For any $k$ large enough satisfying
$$\big|\Omega_{k,0} \big|_{{k,AF}} =
\int_{\Omega_{k,0}}G_k^4 dA_{k} \leq (\alpha/2)^4 s_1,$$ we must have that
$$\Omega_{k,0} \cap G_{k}^{-1}[0, \alpha/2] \neq \emptyset.$$

We just need to check that the distance between $A_k$ and $B_k$
is bounded from below independently of $k$.
Let
$$d := \dist_{{g_\infty}}
\left(G_{\infty}^{-1}[\alpha,\,+ \infty),\,G_{\infty}^{-1}[0, \alpha/2]\right) > 0.$$
Regarding $B_k$ as its lift to $X_{\infty}$,
we have from Lemma~\ref{greens} (ii),
$$B_k \subset G_{\infty}^{-1}[0, \alpha/2].$$
Therefore, for all $k$ sufficiently large,
$$\dist_{{k,AF}}(A_k, B_k) \geq
\dist_{{\infty,AF}}(A_k, G_{\infty}^{-1}[0,\alpha/2]) \geq (\alpha/2)^2 d.$$
\end{proof}

\begin{proof}[{\bf Proof of Lemma \ref{greens_small}}]

For any $\delta > 0$, set
$$\eta_{k,\delta} := \min\left\{1, \frac{(G_k - \delta)_+}{\delta}\right\},$$
where $(G_k - \delta)_+$ denotes the positive part of $(G_k - \delta)$.
First note that $\partial X_k $ is minimal
and its area goes to zero when $k$ goes to infinity.
Using a similar arguments to the proof of Lemma~\ref{essential},
we obtain that, for all $k$ sufficiently large,
$$ \partial X_k \subset G_k^{-1}[0,\,\delta].$$

It suffices to check that, for any $\varepsilon > 0$,
we can find a small $\delta > 0$ so that
$$\int_{X_k}|\nabla \eta_{k,\delta}|^2\,dV_{k,AF} \leq \varepsilon \qquad \mbox{and}
\qquad \vol_{k,AF}\Big(\Big\{x \in X_k\,|\,\nabla \eta_{k,\delta}\neq 0\}\Big) \leq \varepsilon$$
for all $k$ sufficiently large.

For any open set $O$ containing $p_{\infty}$ in $N_{\infty}$,
Lemma~\ref{greens} (iii) and (v) imply that
$$\int_{N_{\infty}-O}G^6_{\infty}\,dV_{\infty} < \infty.$$
Combining this with Lemma~\ref{greens} (iii), we can choose small $\delta > 0$ so that
$$\int_{G^{-1}_{\infty}[\delta,\,2\delta]}G^6_{\infty}\,dV_{\infty}
\leq \frac{\varepsilon}{2}.$$
The uniform approximation property of Lemma~\ref{greens} (ii) guarantees that
$$\vol_{k,AF}\Big(\{x \in X_k\,|\,\nabla \eta_{k,\delta}\neq 0\}\Big) \leq
\int_{G^{-1}_{k}[\delta,\,2\delta]}G^6_{k}\,dV_{k} < \varepsilon$$
for all $k$ large enough.
Likewise, Lemma~\ref{greens} (iv) implies that
\begin{align*}
\int_{X_k}|\nabla \eta_{k,\delta}|_{k,AF}^2\,dV_{k,AF} & \leq
\int_{G^{-1}_{k}[\delta,\,2\delta]}G_k^2\delta^{-2}|\nabla G_k|^2\,dV_k\\
&  \leq
4L^2 \int_{G^{-1}_{k}[\delta,\,2\delta]}G^2_{k}\,dV_k.
\end{align*}
Arguing as before, we can find a small $\delta > 0$ so that
$$4L^2\int_{G^{-1}_{k}[\delta,\,2\delta]}G^2_{k}\,dV_k
< \varepsilon.$$
for all $k$ sufficiently large.
This completes the proof.
\end{proof}

\section{Proof of main theorem: general case}\label{general}

In this section, we complete the proof of Theorem~2.1. {Let $L$ be a closed $3$-manifold for which each prime
factor contains either a non-separating $S^2$ or a two-sided $\RP^2$}. Then, according to \cite[Theorem~12.7 and
its Proof]{Hempel}, $L$ can be expressed as
$$L := \# \ell_1 (S^2 \times S^1) \# \ell_2(S^2 \tilde{\times} S^1) \# Q_1 \# \cdots \# Q_{\ell_3}, $$
where $\ell_1 + \ell_2 + \ell_3 \geq 1$ and each $Q_i$ is a prime closed $3$-manifold with a two-sided $\RP^2$
which is homotopy equivalent to $\RP^2 \times S^1$. Furthermore, set
$$N := M \# L,$$
where $M$ is a closed $3$-manifold satisfying the conditions (ii) and (iii) in Theorem~2.1. For the proof of the
assertion (1) in Theorem~2.1, we assume that $M$ is not diffeomorphic to $S^3$ and $\ell_1 + \ell_2 \geq 1$. For
the proof of the assertion (2) in Theorem~2.1, { we} assume that $\ell_3 \geq 1$. In order to prove Theorem
\ref{main} it suffices to show that $$Y(N, [g]) < Y_2$$ for any unit-volume Yamabe metric $g$ on $N$ with $Y(N,
[g]) > 0$.

We now describe a suitable infinite Riemannian covering of $(N,g)$.
If $L$ has no two-sided $\RP^2$ (i.e., $\ell_3 = 0$),
define $L_{\infty}$ to be the (infinite) universal covering of $L$.
In this case, $\ell_1 + \ell_2 \geq 1$ and
$M$ satisfies the conditions (i), (ii) and (iii) in Theorem~2.1.
If $L$ has a two-sided $\RP^2$ (i.e., $\ell_3 \geq 1$),
then there exists a normal infinite covering
$L_{\infty} \rightarrow L$ having a two-sided $\RP^2$ in its fundamental domain and such that
every embedded $2$-sphere and every two-sided $\RP^2$
in $L_{\infty}$ separate $L_{\infty}$.
Finally, let $(N_{\infty},g_{\infty})$ be the normal infinite Riemannian covering of $(N,g)$
associated to the normal infinite covering $L_{\infty}$ of $L$ which is topologically
$$N_{\infty} = \#_1^{\infty}M \# L_{\infty}.$$
Before stating the relevant properties of $N_{\infty}$,
we need to introduce some notation.

Fix $p_{\infty} \in N_{\infty}$ and consider a sequence of finite Riemannian coverings $\{(N_k, g_k)\}_{k \geq
1}$ of $(N,g)$ satisfying $(N_1, g_1) = (N, g)$ and the following:
\begin{enumerate}
\item[(i)] $N_{\infty}$ is an infinite covering of each $N_k$.

\item[(ii)] $N_{k + 1}$ is a finite covering of $N_k$ for $k \geq 1$.

\item[(iii)] For $k \geq 2$, there exists a fundamental domain of $N_k$ containing $p_{\infty}$ such that its
closure $\widehat{N_{k}}\ (\subset N_{\infty})$ satisfies
$$\{x \in N_{\infty}\ |\ \mbox{dist}_{g_{\infty}}(x,p_{\infty}) \leq k\}
\subsetneq \widehat{N_{k}}.$$
\end{enumerate}

We can assume, without loss of generality,
that $\partial \widehat{N_k}$ is smooth
and that $\partial \widehat{N_{k+1}}, \partial \widehat{N_k}$
have the same projection on $N_k$ for all $k\geq 2$.
Denote by $\widehat{N}\ (\subset N_{\infty})$
the closure of a fundamental domain of $N$ containing $p_{\infty}$
such that $\partial \widehat{N_2}$ and $\partial \widehat{N}$
have the same projection on $N$.

The important properties of $N_{\infty}$ are
that every embedded $2$-sphere and every two-sided $\RP^2$
in $N_{\infty}$ separate $N_{\infty}$,
and that  $\widehat{N}-\{p_{\infty}\}$ is not diffeomorphic to $\R^3$ minus a finite number of balls.
This last assertion follows from the fact that either $M$ satisfies the conditions (i)-(iii) in Theorem~2.1
or the covering $N_{\infty} \rightarrow N$ has a fundamental domain
with a two-sided $\RP^2$.

Let $$\mathcal{G} := \pi_1(N)/\pi_1(N_{\infty})\quad
\mbox{and}\quad\mathcal{G}_k := \pi_1(N_k)/\pi_1(N_{\infty})\subset \mathcal{G}$$
denote the groups
of deck transformations for the normal covering $N_{\infty} \rightarrow N$
and  for each of the normal coverings $N_{\infty} \rightarrow N_k$, respectively.
Here, we identify $\pi_1(N_{\infty})$ with
its projections to $\pi_1(N)$ and $\pi_1(N_k)$.
Set
$$W_n :=
\{\gamma \in \mathcal{G}\ |
\ \mbox{dist}_{g_{\infty}}(\gamma(p_{\infty}),p_{\infty})\geq n\}.$$
Each $\mathcal{G}_k$ has the property that
 $$\mathcal{G}_k - \{e\} \subset W_{q(k)},$$
where $\{q(k)\}$ is a sequence of positive integers
going to infinity with $q(k) \leq q(k+1)$ for all $k \geq 1$.

We denote the projection of $p_{\infty}$ to each $N_k$
also by $p_{\infty}$.
Like in Section~4,
let $G_k$ be the normalized Green's function on $N_k$
for the conformal Laplacian $L_{g_k}$ with the pole at $p_{\infty}$.
Consider the following scalar-flat, \as metrics
$$g_{k,AF} := G_k^4g_k \quad \mbox{on}\quad N_k^* := N_k - \{p_{\infty}\}.$$

Using Aubin's Lemma~3.6,
$$Y(N, [g]) < Y(N_k, [g_{k}]) < Y(N_{k+1}, [g_{k+1}])$$
for all $k \geq 2$.
Therefore, it suffices to prove the following theorem.

\begin{thm}\label{main4}
 $$ \lim_{k \to \infty} Y(N_k^*,[g_{k,AF}]) \leq Y_2.$$
\end{thm}

\begin{proof}
The strategy for the proof is the same as in Sections~4 and 5.

Consider the normalized minimal positive Green's function $G_{\infty}$
on $N_{\infty}$ for $L_{g_{\infty}}$ with the pole at $p_{\infty}$.
This function is defined on $N^*_{\infty} := N_{\infty}-\{p_{\infty}\}$.
Note that Lemma~5.1 still holds,
with some modifications, in this new setting.
For example, Lemma~5.1 (i) becomes
\begin{equation}\tag{6}\label{p1}
G_k=\sum_{\gamma \in \mathcal{G}_k} G_{\infty}\circ \gamma \quad\mbox{on }N_{\infty}
\end{equation}
for all $k\geq 1$, and Lemma~5.1 (iii) now says that
\begin{equation}\tag{7}\label{p2}
\lim_{n \to \infty} \sup  \{G_{\infty}(\gamma(x))\,\,|\,\, x \in K, \,\gamma \in W_n \} = 0
\end{equation}
for every compact set $K$ of $N^*_{\infty}$.

The main difference between the case we consider in this section
and the one considered in Section~5
is that the number of ``necks'' in $N_k^*$ increases with $k$.
Nevertheless, we show that
these ``necks'' are very small when $k$ is large enough.

\begin{lemm}
For any $k$, there exists a family $\Gamma_k$ of disjoint embedded surfaces in $N_k^*$
such that each of its connected components is either a $2$-sphere
or a two-sided $\RP^2$ which minimizes area $($with respect to $g_{k,AF})$, and that
$$\lim_{k \to \infty}\big |\Gamma_k \big|_{k,AF}=0.$$
Moreover, every embedded $2$-sphere and
every embedded two-sided $\RP^2$ of $N_k^*-\Gamma_k$
bound a connected region.
\end{lemm}

\begin{proof}
Take a subset $S$ of the connected components of $\partial \widehat{N}$ which projects injectively into $N^*$,
such that every embedded $2$-sphere and every embedded two-sided $\RP^2$ of $N - S$ separate. Note that $N - S$
is connected and that $S$ is the disjoint union of connected components $\{S_j\}_{j=1}^{j_0}$, which are either
embedded $2$-spheres or two-sided real projective planes. From the construction of $N_k$, (replacing $\{q(k)\}$
by another sequence with the same properties if necessary) there exist subsets
$$I_{j,k} \subset W_{q(k)} \subset \mathcal{G}\quad \mbox{for}\ \ j=1,\cdots,j_0,$$
such that
$$(I_{1,k}(S_1) \sqcup \cdots \sqcup I_{j_0,k}(S_{j_0})) \subset \partial \widehat{N_k},$$
and every embedded $2$-sphere and every two-sided $\RP^2$ of
$$N_k-(I_{1,k}(S_1) \sqcup \cdots \sqcup I_{j_0,k}(S_{j_0}))$$
separate,
where
$$I_{j,k}(S_j) := \{\gamma(S_j)\,\,|\,\, \gamma \in I_{j,k}\}
\quad \mbox{for}\ \ j=1,\cdots,j_0.$$
Note also that $N_k-(I_{1,k}(S_1) \sqcup \cdots \sqcup I_{j_0,k}(S_{j_0}))$ is connected.
We now check that
$$\lim_{k \to \infty}\big |I_{1,k}(S_1)\sqcup \cdots \sqcup I_{j_0,k}(S_{j_0}) \big |_{k,AF}=0.$$
From (6), we have for $j=1,\cdots, j_0$,
$$\sum_{\gamma \in \mathcal{G}}\int_{S_{j}} G_{\infty}\circ \gamma \,dA_{\infty}
= \int_{S_{j}}G_1\, dA_1 < \infty.$$ Hence,
\begin{align*}
\lim _{k \to \infty}\int_{I_{j,k}(S_j)}G_k\,dA_k & =  \lim_{k \to \infty}\sum_{\gamma \in \mathcal{G}_k}
\int_{I_{j,k}(S_j)} G_{\infty}\circ \gamma \,dA_{\infty}\\
& \leq  \lim_{k \to \infty }\sum_{\gamma \in W_{q(k)}}
\int_{S_{j}} G_{\infty} \circ \gamma \,dA_{\infty}\\
&= 0,
\end{align*}
provided that $q(k)$ goes to infinity.
Using \eqref{p2}, we have that
$$\lim_{k \to \infty} \sup \{G_{\infty}(x)\,\,|\,\, x \in I_{i,k}(S_i)\}
=  0 \quad \mbox{for}\ \ j=1,\cdots,j_0,$$
and thus this implies
\begin{align*}\lim_{k \to \infty}\big |I_{j,k}(S_j) \big |_{k,AF}
& = \lim_{k \to \infty}\int_{I_{j,k}(S_j)}G_k^4\, dA_{k} \\
& \leq  \lim _{k \to \infty}\int_{I_{j,k}(S_j)}G_k \,dA_k \\
& =  0.
\end{align*}
The proof follows by minimizing area (with respect to $g_{k,AF}$)
in the isotopy class of
$I_{1,k}(S_1) \sqcup \cdots \sqcup I_{j_0,k}(S_{j_0})$.
The existence and regularity theory in \cite{MSY} says that
the area minimizer is, in this case, a disjoint union of $2$-spheres
and real projective planes.
Again, according to \cite[Section~3]{MSY},
the connected components of the area minimizer
can be jointed by arbitrarily thin tubes so that
the resulting surface is isotopic to
$I_{1,k}(S_1) \sqcup \cdots \sqcup I_{j_0,k}(S_{j_0})$.
Hence, the area minimizer has to contain
a minimal surface $\Gamma_k\ (\subset N_k^*)$
consisting of $2$-spheres or two-sided real projective planes such that
every embedded $2$-sphere and every two-sided $\RP^2$
in $N_k^* - \Gamma_k$ separate.
\end{proof}

Set $X_k$ to be the metric completion of
$$N_k^* - \Gamma_k = N_k - (\{p_{\infty}\}\cup \Gamma_k),$$
and denote the smooth metric extension of $g_{k,AF}$ to $X_k$ by the same symbol.
Then, $(X_k,g_{k,AF})$ becomes an \as manifold
with minimal boundary consisting of a disjoint union of $2$-spheres and projective planes.
Note that Theorem~3.4 and Lemma~3.5
imply the monotonicity of the Hawking quasi-local mass
for each solution of the \imc flow on $(X_k, g_{k,AF})$,
provided that its initial surface is the boundary of a connected region.

Since  $\widehat{N}-\{p_{\infty}\}$ is not diffeomorphic to $\R^3$ minus a finite number of balls, there exists
a positive constant $\alpha$ such that $G_{\infty}^{-1}[\alpha, + \infty)$ has { also this} property. The proof
of Lemma 5.4 can then be easily modified to give the following.
\begin{lemm}
The outermost minimal surface of each $(X_k, g_{k,AF})$ has a connected
component $\Sigma(k)$ with area bounded from below independently of $k$.
\end{lemm}

The argument in the proof of Lemma~5.6 is
still applicable to this new setting.
The same arguments as in the proof of Theorem~4.1 in Section~4
completes the proof of Theorem~7.1.
\end{proof}

\newpage

\bibliographystyle{amsbook}

\vspace{20mm}

\end{document}